# Numerology

by Florentin Smarandache
University of New Mexico
Gallup, NM 87301, USA

**Abstract**. A collection of original sequences, open questions, and problems are mentioned below.

**Introduction**.
One presents many Concatenated Sequences, P-Q relationships, Digital Sequences, Magic Squares, Prime Conjectures, k-Divisibility and Strong Divisibility Sequences, Geometric Sequence, Proposed problems.

### I. Concatenated and operation sequences

1) Odd Sequence:

1, 13, 135, 1357, 13579, 1357911, 135791113, 13579111315, 1357911131517, . . .

How many of them are primes?

2) Even Sequence:

2, 24, 246, 2468, 246810, 24681012, 2468101214, 2468101214, 246810121416, . . .

Conjecture: No number in this sequence is an even power.

3) Prime Sequence:

2, 23, 235, 2357, 235711, 23571113, 2357111317, 235711131719, . . .

How many of them are prime?

Conjecture: A finite number.

4) S-sequence:

General definition: Let $S = \{ s_1, s_2, s_3, \ldots, s_n, \ldots \}$ be an infinite sequence of integers.

Then the corresponding S-sequence is $\{ s_1, s_1 s_2, \ldots, s_1 s_2 \ldots s_n, \ldots \}$ where the numbers are concatenated together.

Question 1: How many terms of the S-sequence are found in the original set S?
Question 2: How many terms of the S-sequence satisfy the properties of other given sequences?

For example, the odd sequence above is built from the set S = { 1, 3, 5, 7, 9, ... } and every element of the S-sequence is found in S. The even sequence is built from the set S = { 2, 4, 6, 8, 10, ... } and every element of the corresponding S-sequence is also in S. However, the question is much harder for the prime sequence.

Study the case when S is the Fibonacci numbers { 1, 1, 2, 3, 5, 8, 13, 21, ... }. The corresponding F-sequence is then { 1, 11, 112, 1123, 11235, 112358, 11235813, ... }. In particular, how many primes are in the F-sequence?

5) Uniform sequences:

General definition: Let $n \neq 0$ be an integer and $d_1, d_2, \ldots, d_r$ distinct digits in a base $B > r$. Then, multiples of n, written using only the complete set of digits $d_1, d_2, \ldots, d_r$ in base B, increasingly ordered, is called the uniform sequence.

Some particular examples involve one digit only.

a) Multiples of 7 written in base 10 using only the digit 1.

111111, 111111111111, 111111111111111111, 111111111111111111111111, ...

b) Multiples of 7 written in base 10 using only the digit 2.

222222, 222222222222, 222222222222222222, 222222222222222222222222, ...

c) Multiples of 79365 written in base 10 using only the digit 5.

55555, 555555555555, 555555555555555555, 555555555555555555555555, ...

In many cases, the uniform sequence is empty.

d) It is not possible to create multiples of 79365 in base 10 using only the digit 6.

Remark: If there exists at least one such multiple of n written with the digits $d_1, d_2, \ldots, d_r$ in base B, then there exists an infinite number of multiples of n. If m is the initial multiple, then they all have the form, m, mm, mmm, ... .

With a computer program it is easy to select all multiples of a given number written with a set of digits, up to a maximum number of digits.

Exercise: Find the general term expression for multiples of 7 using only the digits { 1, 3, 5 } in base 10.

6) Operation Sequences:

General definition: Let E be an ordered set of elements, $E = \{ e_1, e_2, \ldots \}$ and $\Theta$ a set of binary operations well-defined on E. Then

$a_1 \in \{e_1, e_2, \ldots\}$

$a_{n+1} = \min\{e_1 \theta_1 e_2 \theta_2 \ldots \theta_n e_{n+1}\} > a_n$, for $n \geq 1$.

where all $\theta_i$ are operations belonging to $\Theta$.

Some examples:

a) When E is the set of natural numbers and $\Theta = \{+, -, *, /\}$, the four standard arithmetic operations. Then

$a_1 = 1$

$a_{n+1} = \min\{1 \theta_1 2 \theta_2 \ldots \theta_n (n+1)\} > a_n$, for $n \geq 1$.

where $\theta_i \in \Theta$.

Questions:

a) Given N as the set of numbers and $\Theta = \{+,-,*,/\}$ as the set of operations, is there a general formula for this sequence?

b) If the finite sequence is defined with the finite set of numbers $\{1, 2, 3, \ldots, 99\}$ and the set of operations the same as above, where

$a_1 = 1$

$a_{n+1} = \min\{1 \theta_1 2 \theta_2 \ldots \theta_{98} 99\} > a_n$, for $n \geq 1$.

Same questions as in (a).

c) Let N be the set of numbers and $\Theta = \{+, -, *, /, **, (\sqrt{})\}$, where x**y is x to the power y and $x(\sqrt{})y$ is the xth root of y. Define the sequence by

$a_1 = 1$

$a_{n+1} = \min\{1 \theta_1 2 \theta_2 \ldots \theta_n (n+1)\} > a_n$, for $n \geq 1$.

The same questions can be asked, although they are harder and perhaps more interesting.

d) Using the same set of operations, the algebraic operation finite sequence can be defined:

$a_1 = 1$

$$a_{n+1} = \min\{1 \theta_1 2 \theta_2 \ldots \theta_{98} 99\} > a_n, \text{ for } n \geq 1.$$

And pose the same questions as in (b).

More generally, the binary operations can be replaced by $k_i$-ary operations, where all $k_i$ are integers.

$$a_1 \varepsilon \{e_1, e_2, \ldots\}$$

$$a_{n+1} = \min\{1 \theta_1 2 \theta_1 \ldots \theta_1 k_1 \theta_2 (k_1+1) \theta_2 \ldots \theta_2 (k_1 + k_2 - 1) \ldots$$
$$(n+2-k_r) \theta_r \ldots \theta_r (n+1)\} > a_n$$

where $n \geq 1$.

Where each $\theta_i$ is a $k_i$-ary relation and $k_1 + (k_2 - 1) + \ldots + (k_r - 1) = n + 1$.

Note that the last element of the $k_i$ relation is the first element of the $k_{i+1}$ relation.

Remark: The questions are much easier when $\Theta = \{+,-\}$. Study the operation type sequences in this easier case.

e)  Operation sequences at random:

Same definitions and questions as the previous sequences, except that the minimum condition is removed.

$$a_{n+1} = \{e_1 \theta_1 e_2 \theta_2 \ldots \theta_n e_{n+1}\} > a_n, \text{ for } n \geq 1.$$

Therefore, $a_{n+1}$ will be chosen at random, with the only restriction being that it be greater than $a_n$.

Study these sequences using a computer program with a random number generator to choose $a_{n+1}$.

Reference:
F.  Smarandache, "Properties of the Numbers", University of Craiova Archives, 1975. [ Also see the Arizona State University Special Collections, Tempe, Arizona, USA ].

## II.  P - Q Relationships and Sequences

Let $A = \{a_n\}$, $n \geq 1$ be a sequence of numbers and p, q integers $\geq 1$. We say that the terms

$$a_{k+1}, a_{k+2}, \ldots, a_{k+p}, a_{k+p+1}, a_{k+p+2}, \ldots, a_{k+p+q}$$

satisfy a p-q relationship if

$$a_{k+1} \lozenge a_{k+2} \lozenge \ldots \lozenge a_{k+p} = a_{k+p+1} \lozenge a_{k+p+2} \lozenge \ldots \lozenge a_{k+p+q}$$

where $\lozenge$ may be any arithmetic or analytic operation, although it is generally a binary relation on A.

If this relationship is satisfied for any $k \geq 1$, then $\{a_n\}$, $n \geq 1$ is said to be a p - q - $\lozenge$ sequence. For operations such as addition, where $\lozenge = +$, the sequence is called a p - q -additive sequence.

As a specific case, we can easily see that the Fibonacci/Lucas sequence ($a_n + a_{n+1} = a_{n+2}$, for $n \geq 1$), is a 3 - 1 - additive sequence.

Definition: Given any integer $n \geq 1$, the value of the Smarandache function $S(n)$ is the smallest integer m such that n divides m!.

If we consider the sequence of numbers that are the values of the Smarandache function for the integers $n \geq 1$,

1, 2, 3, 4, 5, 3, 7, 4, 6, 5, 11, 4, 13, 7, 5, 6, 17, . . .

they can be incorporated into questions involving the p - q - $\lozenge$ relationships.

a)  How many ordered quadruples are there of the form (S(n), S(n+1), S(n+2),S(n+3)) such that

   S(n+1) + S(n+2) = S(n+3) + S(n+4)

which is a 2 - 2 - additive relationship?

 The three quadruples

   S(6) + S(7) = S(8) + S(9) ,   3 + 7 = 4 + 6;
   S(7) + S(8) = S(9) + S(10),  7 + 4 = 6 + 5;
   S(28) + S(29) = S(30) + S(31), 7 + 29 = 5 + 31.

are known. Are there any others? At this time, these are the only known solutions.

b)  How many quadruples satisfy the 2 - 2 - subtractive relationship

   S(n+1) - S(n+2) = S(n+3) - S(n+4)?
The three quadruples

   S(1) - S(2) = S(3) - S(4),  1 - 2 = 3 - 4;

$$S(2) - S(3) = S(4) - S(5), \quad 2 - 3 = 4 - 5;$$
$$S(49) - S(50) = S(51) - S(52), \quad 14 - 10 = 17 - 13$$

are known. Are there any others?

c) How many 6-tuples satisfy the 3 - 3 - additive relationship

$$S(n+1) + S(n+2) + S(n+3) = S(n+4) + S(n+5) + S(n+6)?$$

The only known solution is

$$S(5) + S(6) + S(7) = S(8) + S(9) + S(10), \quad 5 + 3 + 7 = 4 + 6 + 5.$$

Charles Ashbacher has a computer program that calculates the values of the Smarandache function. Therefore, he may be able to find additional solutions to these problems.

More generally, if $f_p$ is a p-ary relation and $g_q$ a q-ary relation, both defined on the set $\{a_1, a_2, a_3, \dots\}$, then

$$a_{i_1}, a_{i_2}, \dots, a_{i_p}, a_{j_1}, a_{j_2}, \dots, a_{j_q}$$

satisfies a $f_p - g_q$ relationship if

$$f_p(a_{i_1}, a_{i_2}, \dots, a_{i_p}) = g_q(a_{j_1}, a_{j_2}, \dots, a_{j_q}).$$

If this relationship holds for all terms of the sequence, then $\{a_n\}$, $n \geq 1$ is called a $f_p - g_q$ sequence.

Study some $f_p - g_q$ relationships for well-known sequences, such as the perfect numbers, Ulam numbers, adundant numbers, Catalan numbers and Cullen numbers. For example, a 2 - 2 - additive, subtractive or multiplicative relationship.

If $f_p$ is a p-ary relationship on $\{a_1, a_2, a_3, \dots\}$ and

$$f_p(a_{i_1}, a_{i_2}, \dots, a_{i_p}) = f_p(a_{j_1}, a_{j_2}, \dots, a_{j_p})$$

for all $a_{i_k}, a_{j_k}$ where $k = 1, 2, 3, \dots, p$ and for all $p \geq 1$, then

$\{a_n\}$, $n \geq 1$ is called a perfect f - sequence.

If not all p-plets $(a_{i_1}, a_{i_2}, \ldots, a_{i_p})$ and $(a_{j_1}, a_{j_2}, \ldots, a_{j_p})$ satisfy the f relation or the relation is not satisfied for all $p \geq 1$, then $\{a_n\}$, $n \geq 1$ is called a partial perfect f - sequence.

For example, the sequence

1, 1, 0, 2, -1, 1, 1, 3, -2, 0, 0, 2, 1, 1, 3, 5, -4, -2, -1, 1, -1, 1, 3, 0, 2, . . .

is a partial perfect additive-sequence. This sequence has the property that

$$\sum_{i=1}^{p} a_i = \sum_{j=p+1}^{2p} a_j, \text{ for all } p \geq 1.$$

It is constructed in the following way:

$a_1 = a_2 = 1,$

$a_{2p+1} = a_{p+1} - 1$

$a_{2p+2} = a_{p+1} + 1$

for all $p \geq 1$.

a) Can you, the reader, find a general expression of $a_n$ (as a function of n)? Is it periodic, convergent or bounded?

b) Develop other perfect or partial perfect f - sequences. Think about multiplicative sequences of this type.

References:

[1] Sloane, N. J. A, Plouffe, Simon, " The Encyclopedia of Integer Sequences", Academic Press, San Diego, New York, Boston, Sydney, Tokyo, Toronto, 1995/M0453.
[2] Smarandache, F., "Properties of the Numbers", 1975, University of Craiova Archives; (See also Arizona State University Special Collections, Tempe, AZ, USA.)

### III. Digital Subsequences

Let $\{a_n\}$ $n \geq 1$ be a sequence defined by a property (or a relationship involving its terms) P. We then screen this sequence, selecting only the terms whose digits also satisfy the property or relationship.

1) The new sequence is then called a P - digital subsequence.

Examples:

a) Square-digital subsequence:

Given the sequence of perfect squares

   0, 1, 4, 9, 16, 25, 36, 49, 64, 81, 100, 121, 144, ...

only those terms whose digits are all perfect squares { 0, 1, 4, 9 } are chosen. The first few terms are

   0, 1, 4, 9, 49, 100, 144, 400, 441.

Disregarding squares of the form N00 ... 0, where N is also a perfect square, how many numbers belong to this subsequence?

b)   Given the sequence of perfect cubes,

 0, 1, 8, 27, 64, 125, . . .

only those terms whose digits are all perfect cubes { 0, 1, 8 } are chosen. The first few terms are

 0, 1, 8, 1000, 8000.

Disregarding cubes of the form N00 . . . 0, where N is also a perfect cube, how many numbers belong to this subsequence?

c)   Prime-digital subsequence:

Given the sequence of prime numbers,

2, 3, 5, 7, 11, 13, 17, 19, 23, . . .

Only those primes where all digits are prime numbers are chosen. The first few elements of this subsequence are

   2, 3, 5, 7, 23, 29, . . .

Conjecture: This subsequence is infinite.

 In the same vein, elements of a sequence can be chosen if groups of digits, except the complete number, satisfy a property (or relationship) P.  The subsequence is then called a P-partial-digital subsequence.
 Examples:

a) Square-partial-digital subsequence:

49, 100, 144, 169, 361, 400, 441, . . .

In other words, perfect squares whose digits can be partitioned into two or more groups that are perfect squares.

For example 169 can be partitioned into 16 and 9.

Disregarding the square numbers of the form $N00..0$, where N is a perfect square, how many other numbers belong to this sequence?

b) Cube-partial-digital subsequence:

1000, 8000, 10648, 27000, . . .

i.e. all perfect cubes where the digits can be partitioned into two or more groups that are perfect cubes. For example 10648 can be partitioned into 1, 0, 64 and 8.

Disregarding the cube numbers of the form $\overline{N000...00}$, where N is a perfect cube, how many other numbers belong to this sequence?

c) Prime-partial-digital subsequence:

23, 37, 53, 73, 113, 137, 173, 193, 197, . . .

i.e. all prime numbers where the digits can be partitioned into two or more groups of digits that are prime numbers. For example, 113 can be partitioned into 11 and 3.

Conjecture: This subset of the prime numbers is infinite.

d) Lucas-partial-digital subsequence:

Definition: A number is a Lucas number if it is a member of the sequence

$L(0) = 2, L(1) = 1$ and $L(n+2) = L(n+1) + L(n)$ for $n \geq 1$.

The first few elements of this sequence are

2, 1, 3, 4, 7, 11, 18, 29, 47, 76, 123, 199, . . .

A number is an element of the Lucas-partial-digital subsequence if it is a Lucas number and the digits can be partitioned into three groups such that the third group, moving left to right, is the sum of the first two groups. For example, 123 satisfies all of these properties.

Is 123 the only Lucas number that satisfies the properties of this partition?

Study some P-partial-digital subsequences using the sequences of numbers

i) Fibonacci numbers. A search was conducted looking for Fibonacci numbers that satisfy the properties of such a partition, but none were found. Are there any such numbers?
ii) Smith numbers, Eulerian numbers, Bernouli numbers, Mock theta numbers and Smarandache type sequences are other candidate sequences.

Remark: Some sequences may not be partitionable in this manner.

If a sequence $\{a_n\}$, $n \geq 1$ is defined by $a_n = f(n)$, a function of n, then an f-digital sequence is obtained by screening the sequence and selecting only those numbers that can be partitioned into two groups of digits $g_1$ and $g_2$ such that $g_2 = f(g_1)$.

Examples:

a) If $a_n = 2n$, $n \geq 1$, then the even-digital subsequence is

12, 24, 36, 48, 510, 612, 714, 816, 918, 1020, . . .

where 714 can be partitioned into 7 and 14 in that order and

b) Lucky-digital subsequence:

Definition: Given the set of natural numbers

1, 2, 3, 4, 5, 6, 7, 8, 9, 10, 11, 12, 13, 14, 15, . . .

First strike out every other number, leaving

1, 3, 5, 7, 9, 11, 13, 15, 17, 19, 21, . . .

Then, strike out every third number in the remaining list, every fourth number in what remains after that, every fifth number remaining after that and so on. The set of numbers that remains after this infinite sequence is performed are the Lucky numbers.

1, 3, 7, 9, 13, 15, 21, 25, 31, 33, 37, 43, 49, 51, 63, . . .

A number is said to be a member of the lucky-digital subsequence if the digits can be partitioned into two number mn in that order such that $L_m = n$.

37 and 49 are both elements of this sequence. How many others are there?

Study this type of sequence for other well-known sequences.

Reference:
F. Smarandache, "Properties of the Numbers", University of Craiova Archives, 1975. [See also the Arizona State Special Collections, Tempe, AZ., USA].

## IV. Magic Squares

For n ≥ 2, let A be a set of n elements and l an n-ary relation defined on A. As a generalization of the XVIth-XVIIth century magic squares, we present the magic square of order n. This is a square array of elements of A arranged so that l applied to all rows and columns yields the same result.

If A is an arithmetic progression and l addition, then many such magic squares are known. The following appeared in Durer's 1514 engraving, "Melancholia"

```
16   3   2  13
 5  10  11   8
 9   6   7  12
 4  15  14   1
```

Questions:

1) Can you find a magic square of order at least three or four where A is a set of prime numbers and l is addition?

2) Same question when A is a set of square, cube or other special numbers such as the Fibonacci, Lucas, triangular or Smarandache quotients. Given any m, the Smarandache Quotient q(m) is the smallest number k such that mk is a factorial.

A similar definition for the magic cube of order n, where the elements of A are arranged in the form of a cube of length n.

3) Study questions similar to those above for the cube. An interesting law may be

$$l(a_1, a_2, \ldots, a_n) = a_1 + a_2 - a_3 + a_4 - a_5 + \ldots$$

Reference:

F. Smarandache, "Properties of Numbers", University of Craiova Archives, 1975. [ See also the Arizona State University Special Collections, Tempe, AZ., USA].

### V. Prime Conjecture

Any odd number can be expressed as a sum of two primes minus a third prime, not including the trivial solution p = p + q - q.

For example,

```
1 = 3 + 5 - 7 = 5 + 7 - 11 = 7 + 11 - 17 = 11 + 13 - 23 = . . .
3 = 5 + 11 - 13 = 7 + 19 - 23 = 17 + 23 - 37 =  . . .
5 = 3 + 13 - 11 = . . .
7 = 11 + 13 - 17 = . . .
9 = 5 + 7 - 3 = . . .
11 = 7 + 17 - 13 = . . .
```

a) Is this conjecture equivalent to Goldbach's conjecture? The conjecture is that any odd prime ≥ 9 can be expressed as a sum of three primes. This was solved by Vinogradov in 1937 for any odd number greater than $3^{3^{15}}$

b) The number of times each odd number can be expressed as a sum of two primes minus a third prime are called prime conjecture numbers. None of them is known!

c) Write a computer program to check this conjecture for as many positive numbers as possible.

There are infinitely many numbers that cannot be expressed as the absolute difference between a cube and a square. These are called bad numbers(!)

For example, we have conjectured[1] that 5, 6, 7, 10, 13, and 14 are bad numbers. However, 1, 2, 3, 4, 8, 9, 11, 12, and 15 are not as

$1 = |2^3 - 3^2|$   $2 = |3^3 - 5^2|$   $3 = |1^3 - 2^2|$   $4 = |5^3 - 11^2|$   $8 = |1^3 - 3^2|$

$9 = |6^3 - 15^2|$   $11 = |3^3 - 4^2|$   $12 = |13^3 - 47^2|$   $15 = |4^3 - 7^2|$

a) Write a computer program to determine as many bad numbers as possible. Find an ordered array of a's such that
$$a = |x^3 - y^2|, \text{ for x and y integers} \geq 1.$$

Reference:

F. Smarandache, "Properties of Numbers", University of Craiova Archives, 1975. [ See also the Arizona State University Special Collection, Tempe, AZ. , USA ].

### VI. K-Divisibility and Strong Divisibility Sequences

A sequence of rational integers g is called a **divisibility sequence** if and only if

$n \mid m \Rightarrow g(n) \mid g(m)$

for all positive integers n,m. [See [3] and [4]]

Also, g is called a **strong divisibility sequence** if and only if

$(g(n),g(m)) = g((n,m))$

for all positive integers ,m. [See [1], [2], [3], [4] and [5]]

Of course, it is easy to show that the results of the Smarandache function S(n) is neither a divisibility or a strong divisibility sequence because

$4 \mid 20$ but $S(4) = 4$ does not divide $5 = S(20)$.

a) However, is there an infinite subsequence of integers $M = \{m_1, m_2, \ldots\}$ such that S is a divisibility sequence on M?

b) If $P = \{p_1, p_2, \ldots\}$ is the set of prime numbers, the S is not a strong divisibility sequence on P, because for $i \neq j$ we have

$$(S(p_i), S(p_j)) = (p_i, p_j) = 1 \neq 0 = S(1) = S((p_i, p_j)).$$

And the same question can be asked about P as was asked in part (a).

We introduce the following two notions, which are generalizations of a divisibility sequence and a strong divisibility sequence respectively.

1) A k-divisibility sequence, where $k \geq 1$ is an integer, is defined in the following way:

If $n \mid m \Rightarrow g(n) \mid g(m) \Rightarrow g(g(n)) \mid g(g(m)) \Rightarrow \ldots \Rightarrow \underbrace{g(\ldots(g(n))\ldots)}_{k \text{ times}} \mid \underbrace{g(\ldots g(m))\ldots)}_{k \text{ times}}$

for all positive integers n, m.

For example, $g(n) = n!$ is a k-divisibility sequence.

2) A k-strong divisibility sequence, where $k \geq 1$ is an integer, is defined in the following way:

If $(g(n_1), g(n_2), \ldots, g(n_k)) = g((n_1, n_2, \ldots, n_k))$ for all positive integers $n_1, n_2, \ldots, n_k$.

For example, $g(n) = 2n$ is a k-strong divisibility sequence, because

$$(2n_1, 2n_2, \ldots, 2n_k) = 2 * (n_1, n_2, \ldots, n_k) = g((n_1, n_2, \ldots, n_k)).$$

Remarks: If g is a divisibility sequence and we apply its definition k-times, we get that g is a k-divisibility sequence for any $k \geq 1$. The converse is also true.

If g is a k-strong divisibility sequence, $k \geq 2$, then g is a strong divisibility sequence. This can be seen by taking the definition of a k-strong divisibility sequence and replacing $n_1$ by n and all $n_2, \ldots, n_k$ by m to obtain $(g(n), g(m), \ldots, g(m)) = g((n, m, \ldots, m))$ or $(g(n), g(m)) = g((n, m))$.

The converse is also true, as

$$(n_1, n_2, \ldots, n_k) = ((\ldots((n_1, n_2), n_3), \ldots), n_k).$$

Therefore, we found that:

a) The divisibility sequence notion is equivalent to a k-divisibility sequence, or a generalization of a notion is equivalent to itself.

Is there any paradox or dilemma?

b) The strong divisibility sequence notion is equivalent to the k-strong divisibility sequence notion. As before, a generalization of a notion is equivalent to itself.

Again, is there any paradox or dilemma?

References:

[1] Kimberling, C. "Strong Divisibility Sequences With Nonzero Initial Term.", **The Fibonacci Quarterly**, Vol. 16 (1978): pp. 541-544.

[2] Kimberling, C. "Strong Divisibility Sequences and Some Conjectures." ,**The Fibonacci Quarterly**, Vol. 17 (1979): pp. 13-17.

[3] Ward, M. "Note on Divisibility Sequences.", **Bulletin of the American Mathematical Society**, Vol. 38(1937): pp. 725-732.

[4] Ward, M. "A Note on Divisibility Sequences.", **Bulletin of the American Mathematical Society**, Vol. 45(1939): pp. 334-336.

## VII. Geometric Conjecture

a) Let M be an interior point in an $A_1 A_2 \ldots A_n$ convex polygon and $P_i$ the projection of M on $A_i A_{i+1}$ for i = 1, 2, 3, . . . n.

Then,

$$\sum_{i=1}^{n} \overline{MA_i} \geq c \sum_{i=1}^{n} \overline{MP_i}$$

where c is a constant to be found.

For n = 3, it was conjectured by Erdös in 1935 and solved by Mordell in 1937 and Kazarinoff in 1945. In this case c = 2 and the result is called the Erdös-Mordell Theorem.

Question: What happens in 3-space when the polygon is replaced by a polyhedron?

b) More generally: If the projections $P_i$ are considered under a given oriented angle $\alpha \neq 90$ degrees, what happens with the Erdös-Mordell Theorem and the various generalizations?

c) In 3-space, we make the same generalization for a convex polyhedron

$\sum_{n}^{\underline{\quad}} \quad \sum_{m}^{\underline{\quad}}$

$$\sum_{i=1}^{n} MA_i \geq c_1 \sum_{j=1}^{m} MP_j$$

where $P_j$, $1 \leq j \leq m$ are projections of M on all the faces of the polyhedron.

Furthermore,

$$\sum_{i=1}^{n} \overline{MA_i} \geq c_2 \sum_{k=1}^{r} \overline{MT_k}$$

where $T_k$ $1 \leq k \leq m$ are projections of M on all sides of the polyhedron and $c_1$ and $c_2$ are constants to be determined.

[ Kazarinoff conjectured that for the tetrahedron

$$\sum_{i=1}^{4} \overline{MA_i} \geq 2\sqrt{2} \sum_{i=1}^{4} \overline{MP_i}$$

and this is the best possible.

References:

1. P. Erdös, Letter to T. Yau, August, 1995.
2. Alain Bouvier et Michel George, < Dictionnaire des Mathematiques>, Presses Universitaires de France, Paris, 1979, p. 484.

## VIII.  Proposed Problem

Find a natural number n such that the following holds:

1) A $q_1$ th part of it and $a_1$ more are taken away;

2) A $q_2$ th part of the remainder and $a_2$ more are taken away;

. . . . . . . . . . . . . . . . . . . . . . . . . . . . . . . . . . . . . .

S) A $q_s$ th part of the s-1th remainder and $a_s$ more is taken away,

to form the last remainder $r_s$.

**Solution:**

One can work backwards from the last remainder. This is demonstrated for the case s = 3.

I) $\underbrace{|\text{---------}|\text{---------}|\text{---------}|}_{q_1 \text{ ths} \quad a_1 \quad r_1} \quad (a_1 + r_1) * \dfrac{q_1}{q_1 - 1} = N$

II) $\underbrace{|\text{---------}|\text{---------}|\text{---------}|}_{q_2 \text{ ths} \quad a_2 \quad r_2} \quad (a_2 + r_2) \dfrac{q_2}{q_2 - 1} = r_1 \quad \text{(by notation)}$

III) $\underbrace{|\text{---------}|\text{---------}|\text{---------}|}_{q_3 \text{ ths} \quad a_3 \quad r_3} \quad (a_3 + r_3) \dfrac{q_3}{q_3 - 1} = r_2 \quad \text{(by notation)}$

$q_3 - 1$ of $q_3$ ths

Therefore, $N = \left(a_1 + \left(a_2 + (a_3 + r_3)\dfrac{q_3}{q_3 - 1}\right)\dfrac{q_2}{q_2 - 1}\right)\dfrac{q_1}{q_1 - 1}$.

One can generalize and prove by induction that:

$N = \left(a_1 + \left(a_2 + \ldots + \left(a_{s-1} + (a_s + r)\dfrac{q_s}{q_s - 1}\right)\dfrac{q_{s-1}}{q_{s-1} - 1} \ldots \right)\dfrac{q_2}{q_2 - 1}\right)\dfrac{q_1}{q_1 - 1}$.

Reference:

G. Smarandache, " Properties of the Numbers", University of Craiova Archives, 1975. [ Also see the Arizona State University Special Collections, Tempe, Arizona, USA ].